# SMALL CANCELLATION GROUPS AND TRANSLATION NUMBERS

## ILYA KAPOVICH

ABSTRACT. In this paper we prove that C(4)-T(4)-P, C(3)-T(6)-P and C(6)-P small cancellation groups are translation discrete in the strongest possible sense and that in these groups for any $g$ and any $n$ there is an algorithm deciding whether or not the equation $x^n = g$ has a solution. There is also an algorithm for calculating for each $g$ the maximum $n$ such that $g$ is an $n$-th power of some element. We also note that these groups cannot contain isomorphic copies of the group of $p$-adic fractions and so in particular of the group of rational numbers. Besides we show that for $C''(4) - T(4)$ and $C''(3) - T(6)$ groups all translation numbers are rational and have bounded denominators.

## 0. INTRODUCTION.

In [GS1] S.Gersten and H.Short introduced the notion of translation numbers for finitely generated groups. This concept was inspired by some considerations about groups acting on $\mathbb{R}$-trees and about the action of the fundamental group of a closed riemannian manifold of nonpositive curvature on the universal cover of this manifold (see, for example, the book of W.Ballman, M.Gromov and V.Schroeder "Manifolds of nonpositive curvature" [BGrS]).

**Definition.** Let $G$ be a group and let $X$ be a finite generating set of $G$ closed under inversion. Then any element $g$ in $G$ can be expressed as a product $g = x_1 \cdot x_2 \cdot \ldots \cdot x_n$ where $x_i \in X$; we term the minimal such $n$ the X-length of $g$ and denote it $l_X(g)$.

Then for any $g \in G$ we define the translation number of $g$ with respect to $X$ as follows:
$\tau_X(g) = \lim_{n \to \infty} \frac{l_X(g^n)}{n}$. (It is noted in [GS1] that this limit always exists.)

From now on we will assume that all generating sets for all groups in this section are closed under inversions. It turns out (see [GS1] for proofs) that if $G$ is a group and $X$ is a finite generating set then
(a) $\tau_X(g) = \tau_X(g^{-1}) \leq l_X(g)$;
(b) $\tau_X(g) = \tau_X(hgh^{-1})$;
(c) if $g$ is an element of finite order then $\tau_X(g) = 0$;
(d) $\tau_X(g^n) = n\tau_X(g)$ for any integer $n$ and for any $g \in G$;
(e) for any other generating set $Y$ there are positive constants $C_1$ and $C_2$ such that

$$C_1 \tau_Y(g) \leq \tau_X(g) \leq C_2 \tau_Y(g)$$

and

$$\tau_X(g) = 0 \iff \tau_Y(g) = 0$$

for all $g \in G$;
(f) if $gh = hg$ then $\tau_X(gh) \leq \tau_X(g) + \tau_X(h)$.

---

1991 *Mathematics Subject Classification*. Primary 20F06; Secondary 20F10, 20F32.

This research was supported by the Robert E. Gilleece Fellowship at the CUNY Graduate Center

Note that as property (b) indicates, $\tau_X(\cdot)$ is constant on conjugacy classes of elements of $G$. We define the translation number of a conjugacy class $\alpha$ as $\tau_X(\alpha) = \tau_X(g)$ where $g$ is some element of $\alpha$.

Some properties of the set of all translation numbers of a group $G$ do not depend on the choice of generators and provide us with very interesting invariants of the group.

**Definition.** We say that a finitely generated group $G$ is **translation separable** if for some (and therefore for any) finite generating set $X$ of $G$ any element with translation number 0 has finite order.

One can think of a translation separable group $G$ as a group with the property that for any element $g \in G$ of infinite order the map $n \mapsto g^n$ is a quasiisometric embedding of $\mathbb{Z}$ into $G$ for any word metric on $G$ (see [GS1]). Translation separable groups already form an interesting class and have some nice algebraic properties. For example in such a group two different and not mutually inverse powers of an element of infinite order cannot be conjugate (see [GS1, prop.6.7]) and such a group cannot contain a finitely generated nilpotent subgroup unless this subgroup is virtually abelian (see [GS1, prop. 6.9]). The most remarkable result about translation separable groups is the theorem of S.Gersten and H.Short (see [GS1, prop. 6.6] for proof) which states that *biautomatic groups are translation separable.*

The following definition was suggested by G.Conner in [C1].

**Definition.** A group $G$ is said to be **translation discrete** if it is translation separable and for some (and therefore for any) finite generating set $X$ the set $\tau_X(G) \subseteq \mathbb{R}$ has 0 as an isolated point.

G.Conner showed in [C1] that solvable groups of finite virtual cohomological dimension are not translation discrete unless they are virtually abelian. He constructed, also in [C1], an interesting example of a group for which every translation number is an accumulation point of the set of all translation numbers. It is proved in the same paper that discrete cocompact groups of fixed-point-free isometries of proper convex spaces are translation discrete.

It is natural to consider the class of groups having the property that for some finite generating set $X$ the set $\tau_X(G)$ is a discrete subset of real numbers.

The great disadvantage of this property is that, unlike the previous ones, it is not clear whether it is independent of the choice of generators. The reason is that if we have an infinite sequence of different elements $g_n \in G$ such that the numbers $\tau_X(g_n)$ are all different and converge to some positive (!) number then it might happen (at least theoretically) that for some other generating set $Y$ we have $\tau_Y(g_n) = \tau_Y(g_1)$ for all $n$. However there is a situation when this property does not depend on the choice of generators. It is described in the following definition.

**Definition.** We shall say that a group $G$ is **strongly translation discrete** if it is translation separable and for some (and therefore for any) finite generating set $X$ and for any real number $r$ the number of conjugacy classes $\alpha$ in $G$ with $\tau_X(\alpha) \leq r$ is finite.

It is not hard to see that for a strongly translation discrete group $G$ and for any finite generating set $X$ the set $\tau_X(G)$ is discrete. We can note here some good properties of groups which have some kind of translation discreteness. For example, in a translation discrete group $G$ any element $g$ of infinite order has the property that

$$\sup\{n \mid \text{ there is some } h \in G \text{ such that } h^n = g\}$$

is finite (see propery (d) of translation numbers and the definition of being translation discrete).

Therefore such a group $G$ cannot contain subgroups isomorphic to the additive group of rational numbers or the group of $p$-adic fractions $\mathbb{Q}_p = \{k/p^l \mid k \in \mathbb{Z}, l \in \mathbb{N}\}$. For a strongly translation discrete group $G$ the function $f_X(r)$ counting the number of conjugacy classes with $d_X$-translation numbers less or equal to $r$ is invariant (in the same sense as its Dehn function) under the choice of generators and, possibly, might serve as an invariant of the group. It follows from [Gr, 5.2.C] that word hyperbolic groups are strongly translation discrete (see [S] for a proof). Therefore almost all small cancellation groups (such as C(7), C(5)-T(4), C(4)-T(5), C(3)-T(7) groups), which are known to be word hyperbolic, are strongly translation discrete.

We recall that a group presentation $G = <S|R>$ is said to satisfy a C(p)-condition if any $r \in R$ is cyclically reduced and no relator $r \in R$ is a product of less than $p$ pieces. Here a word $w$ in $S \cup S^{-1}$ is called a piece with respect to this presentation if $w$ is an initial segment of two distinct words $y_1$ and $y_2$ such that $y_i$ is a cyclic permutation of $r_i^{\pm 1}$ for some $r_i \in R$. The formulation of other small cancellation conditions is a little more technical and will be given later (see section 1 below). In [GS2] S.Gersten and H.Short investigated another, in some sense, boundary class of small cancellation groups, namely C(4)-T(4)-P, C(6)-P and C(3)-T(6)-P groups. They showed that these groups are biautomatic (and therefore translation separable) but not necessarily word hyperbolic. (For the definitions and properties of automatic, biautomatic groups, regular languages etc. see [ECHLPT].)

Here we prove the following.

**Theorem 0.1.** *Let $G = <S|R>$ be a group presentation satisfying one of C(4)-T(4)-P, C(6)-P or C(3)-T(6)-P small cancellation conditions and let $X$ be any finite generating set for $G$. Then for any $r \in \mathbb{R}$ the set $\{\alpha | \alpha$ is a conjugacy class in $G$ and $\tau_X(\alpha) \leq r\}$ is finite, that is $G$ is strongly translation discrete.*

The following consequences of Theorem 0.1 are of interest.

**Corollary 0.2.** *Let $G$ and $X$ be as in Theorem 0.1. Then $\tau_X(G)$ is a discrete subset of the real line.*

**Corollary 0.3.** *Let $G$ be as in Theorem 0.1. Then for any $g \in G$ of infinite order $\sup\{n | g = x^n$ for some $x \in G\} < \infty$.*

*In particular $G$ does not contain subgroups isomorphic to the group of $p$-adic fractions $\mathbb{Q}_p = \{k/p^l | k \in \mathbb{Z}, l \in \mathbb{N}\}$ (where $p$ is some prime number) and therefore to the group of rational numbers $\mathbb{Q}$.*

We first obtain a quick proof of Theorem 0.1 in section 2 using the theory of CAT(0)-spaces. Then in sections 3 and 4 we give a different proof such that techniques used in it allow us to obtain, in addition to Theorem 0.1, the following results.

**Theorem 0.4.** *Let $G = <S|R>$ be a group presentation satisfying one of C(4)-T(4)-P, C(6)-P or C(3)-T(6)-P small cancellation conditions and $Y = S \cup S^{-1}$. Then for any freely reduced word $w$ in $Y$ and for any natural number $n$ there is an algorithm which determines whether or not the equation $x^n = w$ has a solution in $G$.*

**Theorem 0.5.** *Let $G = <S|R>$ be a group presentation satisfying one of C(4)-T(4)-P, C(6)-P or C(3)-T(6)-P small cancellation conditions and $Y = S \cup S^{-1}$.*

*Then for any freely reduced word $w$ in $Y$ there is an effective algorithm for calculating $\sup\{n | x^n = w$ in $G$ for some $x \in G\}$.*

**Corollary 0.6.** *If $G = <S|R>$ and $Y = S \cup S^{-1}$ are as in Theorem 0.5, then there is an algorithm which, given a pair of words $w_1, w_2$ in $Y$, decides whether $w_1$ is conjugate to a power of $w_2$ in $G$*

*Proof.* Indeed, we first determine $M = \sup\{n | x^n = w_1$ in $G$ for some $x \in G\}$ and then for each $n = -M, -M+1, \ldots, M$ determine whether $w_1$ is conjugate to $w_2^n$ in $G$ (the group $G$ has solvable conjugacy problem by the result of S.Gersten and H.Short [GS]).

We also show that for $C''(4) - T(4)$ and $C''(3) - T(6)$ groups all translation numbers are rational and have bounded denominators. Besides it turns out that for these groups the language of *all geodesic words* is regular.

**Theorem 0.7.** *Let $G = <S|R>$ be a group presentation satisfying $C''(4) - T(4)$ or $C''(3) - T(6)$ small cancellation condition and $X = S \cup S^{-1}$ then for every $g \in G$ the number $\tau_X(g)$ is rational and, moreover, $2\tau_X(g)$ is an integer. Moreover, for every nontrivial element $g \in G$ the element $g^2$ is conjugate to a periodically geodesic element $h$ that is such that $l_X(h^n) = |n| l_X(h)$ for each $n \in \mathbb{Z}$. Also the set of all $d_X$-geodesic words is a regular language.*

The fact, that in a word hyperbolic group all translation numbers are rational and have bounded denominators, follows from a claim of M.Gromov [Gr, 5.2.C], and it was accurately proved by E.Swensen (see [S], Corollary of Theorem 13). E.Swensen also proves in [S] that for a given word hyperbolic group $G$ with a fixed word metric on it there is a constant $N$ such that for every element $g \in G$ of infinite order $g^N$ is conjugate to a periodically geodesic element. Theorem 0.7 provides a similar statement for $C''(4) - T(4)$ or $C'''(3) - T(6)$ groups (for the word metric corresponding to the standard generating set) and implies that in this case the constant $N$ can be chosen to be equal to two.

G.Conner showed that finitely generated nilpotent groups with word metric have rational translation numbers [C2] and gave a remarkable example of a group with irrational translation numbers [C3]. K.Johnsgard [J1] first obtained the description of geodesics for $C''(4) - T(4)$ and $C'''(3) - T(6)$ groups (our Lemma 3.2 and Lemma 4.2), as well as for $C''(6) - T(3)$-groups, and also observed that the set of all geodesics in such groups is a regular language. In [J2] K.Johnsgard also showed that for these groups the language of lexicographically least geodesics is automatic.

Results analogous to Theorem 0.4 and Theorem 0.5 were proved by S.Lipschutz in [L] for a different class of small cancellation groups, namely for $C'(1/6)$ and $C'(1/4) - T(4)$ groups. However these groups are word hyperbolic (see theorem 4.5 and theorem 4.6 of [LS]) unlike the groups we consider in this paper which are not necessarily word hyperbolic (in particular, they can contain subgroups, isomorphic to $\mathbb{Z} \times \mathbb{Z}$). The relation between powers and conjugacy for several other classes of small cancellation groups was investigated by L.Comerford [Com1], [Com2], L.Comerford and B.Truffault [CT], M.Anshel and P.Stebe [AS] and others.

Our results can be applied to many examples of $C(4) - T(4) - P$ groups provided by C.Weinbaum [W] (for instance the groups of prime alternating knots) and to C(3)-T(6)-P presentations of groups acting on Euclidean buildings in the sense of Bruhat-Tits, constructed in [GS2]. Some other examples of C(3)-T(6) presentations were given by M.El-Mosalamy and S.Pride in [MP].

I am grateful to H.Short and G.Baumslag for useful conversations and help in writing this paper.

## 1. Definitions, notations and preliminary facts

**Diagrams over groups and metrics on groups.**

In this section we shall mainly refer to [GS2] where one can find a detailed discussion on the matter. An excellent overview on small cancellation groups and diagrams over groups can be also found in [LS]. Let $G$ be a group and let $X$ be a finite generating set of $G$ closed under inversion. The **word metric** on $G$ associated to $X$ is defined as follows: $d_X(g, h) = l_X(h^{-1}g)$ for any $g, h \in G$. (It is not hard to see that $d_X$ is indeed a metric on $G$.)

For any word $x_1...x_k$ where $x_i \in X$ we denote the corresponding element $x_1 \cdot ... \cdot x_k$ of $G$ by $\overline{x_1...x_k}$. We say that this word is $d_X$-geodesic if $l_X(\overline{x_1...x_k}) = k$.

We say that a word $w$ in the alphabet $X$ is $m$-locally geodesic if any subword $v$ of $w$ of length at most $m$ is $d_X$-geodesic. If $\alpha$ is a conjugacy class in $G$, we define the $X$-length of $\alpha$ as follows:

$$l_X(\alpha) = \min\{l_X(g) | g \in \alpha\}.$$

For any set $S$ we shall denote the free group on $S$ by $F(S)$.

Let $G = <S|R>$ be a group presentation, where $R$ is a set of cyclically reduced words in $F(S)$ closed under taking cyclic permutations and inverses. We fix $X = S \cup S^{-1}$ as a set of semigroup generators for $G$. We shall associate to such presentation a standard 2-complex $K$ with a single 0-cell, one 1-cell for each $s \in S$ and one 2-cell for each relator $r \in R$. A freely reduced word $w$ of $F(S)$ represents an identity element in $G$ if and only if there is a connected simply connected planar 2-complex $D$ and a map $\phi : (D, \partial D) \to (K, K^{(1)})$ such that

(a) vertices of $D$ go to the vertex of $K$, open $i$-cells map homeomorphically to open $i$-cells of $K$ for $i = 1, 2$;

(b) $\phi$ maps the boundary $\partial D$ to the loop representing the word $w$.

We label each oriented edge of $D$ by some letter of $X$ in a natural way.

Such a $D$ we shall call a **singular disc diagram**. If $D$ does not contain a cut vertex (i.e. a vertex whose removal makes $D$ disconnected), we shall call $D$ a **disc diagram**.

A singular disc diagram $D$ is reduced if there are no 2-cells $R_1$ and $R_2$ in $D$ with boundaries containing a common edge $e$ such that the labels on their boundaries, when reading from edge $e$ clockwise on $R_1$ and anticlockwise on $R_2$, are the same. It is easy to see that any nonreduced singular disc diagram can be transformed into a reduced one without changing its boundary label. In this paper all singular disc diagrams are assumed to be reduced.

We equip each singular disc diagram for our presentation with a **piecewise Euclidean structure** by putting each 2-cell with a label of length $n$ to be isometric to the regular euclidean $n$-gon of side one. (It is always possible to do so if all defining relations have length $\geq 3$.)

We shall refer to disc diagrams (i.e. singular disc diagrams without cut vertices) with piecewise Euclidean structure as **PE-disc diagrams**.

For any vertex $v$ on the boundary $\partial D$ of a PE-disc diagram $D$ put $\sigma_v$ to be the sum of corner angles of adjacent to $v$ and define the **turning angle** at $v$ to be $\tau(v) = \pi - \sigma_v$. We extend the notion of turning angle for singular disc diagrams by putting $\tau(v) = -\infty$ for cut vertices (which always lie on the boundary cycle) and using the already defined notion for each disc component.

**Small cancellation conditions.**

Let $G = <S|R>$ be a group presentation as above. As usual we call a non-empty word $u \in F(S)$ a piece with respect to this presentation if there are two different relators $r_1, r_2 \in R$ such that $r_1 \equiv uv_1$ and $r_2 \equiv uv_2$.

We term a word $u$ a $p/q$-relator if there is a relator $r \in R$ such that $u$ is an initial segment of $r$ and $l(u) = (p/q)l(r)$ where $p$ and $q$ are positive integers.

We say that the presentation above satisfies the $C(p)$-condition if no relator is a product of less than $p$ pieces. We say that it satisfies the $T(q)$-condition if for any sequence $r_1, r_2, \ldots, r_k \in R$ such that $3 \leq k < q$ and $r_i \neq r_{i+1}^{-1}$ at least one of the words $r_i r_{i+1}$ is cyclically reduced (all sub-indices are considered mod $k$). The presentation is said to satisfy the condition $P$ if all pieces have length one and no relator is a proper power. And, finally, the presentation is said to be $C''(p)$ if it is $C(p) - P$ and all relators have length $p$.

It was pointed out in [GS2] that if $D$ is a reduced singular disc diagram for a $C(p) - T(q)$ presentation then each 2-cell of $D$ has at least $p$ sides and each vertex from the topological interior of $D$ meets at least $q$ 2-cells. If in addition our presentation satisfies the $C''(p)$ condition then each 2-cell of $D$ has exactly $p$ sides.

## 2. Groups acting on CAT(0)-spaces

**Definition.** A metric space $(X, d)$ is termed **geodesic** if for any points $x, y \in X$ there exists an isometric map $\alpha \colon [0, d(x, y)] \to X$ such that $\alpha(0) = x$ and $\alpha(d(x, y)) = y$. Such an $\alpha$ is called a geodesic segment in $X$ from $x$ to $y$. A metric space $(X, d)$ is called **proper** if for any $x \in X$ and for any $R \geq 0$ the metric ball $B(x, R) = \{y \in X | d(y, x) \leq R\}$ is compact.

A **geodesic triangle** $\Delta$ in $X$ is a triple of points (vertices) $p, q, r \in X$ together with a choice of three geodesic segments (termed sides of $\Delta$), one joining each pair of vertices. A **comparison triangle** $\Delta'$ for $\Delta$ is a triangle $\Delta'$ in the Euclidean plane $\mathbb{E}^2$ with vertices $p', q', r'$ such that $d(p, q) = d(p', q')$, $d(q, r) = d(q', r')$, $d(p, r) = d(p', r')$. Given a side of $\Delta$ and a point $x$ on it, there is a unique point $x'$ on the corresponding side of a corresponding side of $\Delta'$ such that $d(x, e) = d(x', e')$ for each of the endpoints $e$ of the given side to which $x$ belongs. This point $x'$ is called a **comparison point** for $x$.

**Definition.** A geodesic triangle $\Delta$ in a metric space $(X, d)$ is said to satisfy the CAT(0)-inequality if for any comparison triangle $\Delta'$, for every pair of points $x, y$ on the sides of $\Delta$ and every choice of comparison points $x'$ and $y'$

$$d(x, y) \leq d(x', y').$$

A proper geodesic metric space $(X, d)$ is said to be **nonpositively curved** if for any point $p \in X$ there is a neighborhood $U$ of $p$ such that any geodesic triangle in $U$ satisfies the $CAT(0)$-inequality.

**Definition.** A proper geodesic metric space $(X, d)$ is called a **CAT(0)-space** any geodesic triangle in $X$ satisfies the $CAT(0)$-inequality.

For a basic overview of $CAT(0)$-spaces the reader is referred to [B],[Br] and [BH]. $CAT(0)$-spaces are contractible, have convex distance function (see [B] for definitions) and have the property that any two points can be joined by a unique geodesic segment. It is worth mentioning that any complete simply connected riemannian manifold of nonpositive sectional curvature is a CAT(0)-space (see [BGrS]).

**Definition.** Let $(X, d)$ be a $CAT(0)$-space and $\gamma: X \to X$ be an isometry of $X$ onto itself. We define the **translation length** $\tau(\gamma)$ of $\gamma$ as
$$\tau(\gamma) = \inf\{d(x, \gamma(x)) | x \in X\}.$$

The following properties of translation length are of importance.

**Lemma 2.1.** *Let $(X, d)$ be a $CAT(0)$-space. Then*
  (1) *if $\gamma$ is an isometry of $X$ then $\tau(\gamma) = \tau(\gamma^{-1})$;*
  (2) *if $\gamma$ and $\gamma_1$ are isometries of $X$ then $\tau(\gamma) = \tau(\gamma_1^{-1}\gamma\gamma_1)$;*
  (3) *if $\gamma$ is an isometry of $X$, and there is an $x \in X$ such that $\gamma(x) = x$ then $\tau(\gamma) = 0$;*
  (4) *if $\gamma$ is an isometry of $X$, $n > 0$ and $\gamma^n = id_X$ then there is an $x \in X$ such that $\gamma(x) = x$ and $\tau(\gamma) = 0$.*

*Proof.* Statements (1), (2) and (3) are obvious. We will give a sketch of the proof of (4) and refer the reader to the Fixed Point Theorem of [Br] for details. Let $\gamma^n = id_X$. Pick a point $y \in X$. If $\gamma(y) = y$ then the statement is true. So assume $\gamma(y) \neq y$. Clearly the set $Y = \{y, \gamma(y), .., \gamma^{n-1}(y)\}$ is $\gamma$-invariant, i.e. $\gamma(Y) = Y$. Then the convex hull $Y_1$ of $Y$ is also $\gamma$-invariant and so $\gamma(Y_1) = Y_1$. Notice that $Y_1$ is compact since $X$ is proper.

For any compact convex subset $V$ of $X$ there is a unique point $v \in V$ satisfying the following property: there is a number $r \geq 0$ such that $V \subset B(v, r)$ and if $x \in X$, $R \in \mathbb{R}$ are such that $V \subset B(x, R)$ then $r \leq R$. This $v$ is called the center of $V$.

If in our situation $y_1$ is the senter of $Y_1$ then $\gamma(y_1)$ is the center of $\gamma(Y_1)$. Since $\gamma(Y_1) = Y_1$ and the center is unique, $\gamma(y_1) = y_1$.

The following theorem [BH] demonstrates the geometric meaning of a translation length.

**Theorem 2.2.** *Let $(X, d)$ be a $CAT(0)$-space and a group $G$ act on $X$ by isometries. Suppose the action is cocompact, i.e. the factor space $X/G$ is compact, and properly discontinuous, i.e. for any compact subset $K$ of $X$ the set $\{g \in G | g(K) \cap K \neq \emptyset\}$ is finite. Then for any element $g \in G$ the following holds:*
  (1) *the set $MIN(g) = \{x \in X | d(x, g(x)) = \tau(g)\}$ is nonempty;*
  (2) *if $g$ is of infinite order then $\tau(g) > 0$ and there is a biinfinite geodesic $\alpha: \mathbb{R} \to X$ such that for any $r \in \mathbb{R}$, $g(\alpha(r)) = \alpha(r + \tau(g))$, i.e. $g$ restricted to $\alpha$ is a translation by $\tau(g)$;*
  (3) *if $\beta: \mathbb{R} \to X$ is a biinfinite geodesic such that for some $t > 0$*

$$g(\beta(r)) = \beta(r + t)$$

  *for all $r \in \mathbb{R}$ then $t = \tau(g)$.*

*Sketch of the proof.*

(1) Let $g \in G$. First we notice that the set $MIN(g) = \{x \in X | d(x, g(x)) = \tau(g)\}$ is nonempty. Indeed, suppose $(x_i)_{i \in \mathbb{N}}$ is a sequence of points in $X$ such that $\lim_{i \to \infty} d(x_i, g(x_i)) = \tau(g)$. Since $X/G$ is compact, there is a compact closed ball $B = B(a, R)$ such that $X = GB$. So for any $i$ there is an $h_i \in G$ such that $y_i = h_i(x_i) \in B$. Clearly $d(x_i, g(x_i)) = d(h_i(x_i), h_i g(x_i)) = d(y_i, h_i g h_i^{-1}(y_i))$. Since $B$ is compact, we may assume that $\lim_{i\to\infty} y_i = y \in B$. There is some $M > 0$ such that $d(x_i, g(x_i)) < M$ for all $i$. Thus $d(h_i g h_i^{-1}(y_i)), a) < M + R$ for any $i$ and therefore $B(a, M + R) \cap h_i g h_i^{-1} B(a, M + R) \neq \emptyset$ for all $i$. Since the action of $G$ is properly discontinuous and $B(a, M + R)$ is compact, there is an $i_0$ such that for all $i \geq i_0$ $h_i g h_i^{-1} = h_{i_0} g h_{i_0}^{-1}$. Thus $\tau(g) = \lim_{i\to\infty} d(x_i, g(x_i)) = \lim_{i\to\infty} d(y_i, h_i g h_i^{-1}(y_i)) = \lim_{i\to\infty} d(y_i, h_{i_0} g h_{i_0}^{-1}(y_i)) = d(y, h_{i_0} g h_{i_0}^{-1}(y))$. Therefore $\tau(g) = d(h_{i_0}^{-1}(y), g h_{i_0}^{-1}(y))$ and so $MIN(g)$ is nonempty.

Notice that if $g \in G$ is of infinite order then $\tau(g) > 0$. Indeed, if $\tau(g) = 0$ then, since $MIN(g)$ is nonempty, there is $x \in X$ such that $g(x) = x$. Therefore $g^n(x) = x$ for any integer $n$ which contradicts the fact that $G$ acts properly discontinuously.

(2) Let $g \in G$ be of infinite order, $t = \tau(g) > 0$ and $x \in MIN(g)$. Consider a geodesic triangle $\Delta$ in $X$ with vertices $x = x_0, g(x) = x_1, g^2(x) = x_2$ and a comparison triangle $\Delta'$ with vertices $x_0', x_1', x_2'$ in the euclidean plane. Suppose that $d_(x, g^2(x)) < 2t$. Let $y$ be the midpoint of the geodesic segment $[x, g(x)]$ and $z$ be the midpoint of the geodesic segment $[g(x), g^2(x)]$. Clearly $g(y) = z$ because of the uniqueness of geodesics in $X$. By the definition of $\tau(g)$ $d(y, z) \geq t$. On the other hand, if $y'$ and $z'$ are comparison points for $y$ and $z$ in $\Delta'$ then $d(y', z') = (1/2)d(x_0', x_2') < t$- a contradiction. So $d_(x, g^2(x)) = 2t$ and therefore $g(x)$ lies on a geodesic segment $[x, g^2(x)]$. This implies that $g^n(x)$ lies on the geodesic segment $[g^{n-1}(x), g^{n+1}(x)]$ for all integers $n$. It can now be shown that there is a biinfinite geodesic $c$ such that for all integers $n$ $g^n(x)$ lie on $c$. The proof of the last assertion relies on the following fact about $CAT(0)$-spaces.

If $X$ is a $CAT(0)$-space and $[p, q], [y, z]$ are geodesic segments such that $[p, q] \cap [y, z] = [y, q]$ and $d(y, q) > 0$ then $[p, q] \cup [y, z]$ is a geodesic from $p$ to $z$. Indeed, suppose it is not true and $d(p, z) < d(p, q) + d(q, z)$. Consider a geodesic triangle $\Delta$ in $X$ with vertices $p, q, z$ and a comparison triangle $\Delta'$ in the euclidean plane with vertices $p', q', z'$. Recall that by our assumption $d(p, z) < d(p, q) + d(q, z)$ and therefore for any point $u'$ in the interior of the segment $[p', q']$ we have $d(u', z') < d(u', q') + d(q', z')$. Thus for a corresponding point $u$ on the side $[p, q]$ of $\Delta$ we have

$$d(u, z) \leq d(u', z') < d(u', q') + d(q', z') = d(u, q) + d(q, z).$$

However we know that if $u$ is such that $0 < d(u, q) < d(y, q)$ then $d(u, z) = d(u, q) + d(q, z)$ what gives us a contradiction.

(3) Clearly, if $g$ has a $g$-invariant geodesic on which it acts by a nontrivial translation then $g$ has infinite order. Let $t = \tau(g) > 0$ and $\alpha \colon \mathbb{R} \to X$ be a $g$-invariant geodesic constructed as in (2). This means that for any $r \in \mathbb{R}$ $g(\alpha(r)) = \alpha(r + t)$. Suppose that $\beta \colon \mathbb{R} \to X$ is another geodesic and $s > 0$ is such that for any $r \in \mathbb{R}$ $g(\beta(r)) = \beta(r + s)$. Then $t \leq s$ by the definition of a translation length. We have to show that $t = s$.

Suppose it is not true and $t < s$. Put $a = d(\alpha(0), \beta(0))$. Then for any positive integer $n$ $a = d(g^n(\alpha(0)), g^n(\beta(0))) = d(\alpha(nt), \beta(ns))$, $d(\alpha(0), \alpha(nt)) = nt$ and $d(\beta(0), \beta(ns)) = ns$. Therefore by the triangle inequality we have $|ns - nt| = n|s - t| \leq 2a$ for any positive integer $n$. This contradicts our assumption $t < s$.

**Corollary 2.3.** *Suppose the conditions of Theorem 2.2 are satisfied. Then for any $g \in G$ for any $n > 0$, $\tau(g^n) = n\tau(g)$.*

*Proof.* If $g$ is an element of finite order then $\tau(g) = \tau(g^n) = 0$ by Lemma 2.1(4). Assume now that $g$ is an element of infinite order. Put $t = \tau(g)$. Then by Theorem 2.2 $t > 0$ and there is a $g$-invariant geodesic $\alpha \colon \mathbb{R} \to X$ such that for any $r \in \mathbb{R}$ $g(\alpha(r)) = \alpha(r + t)$. Therefore $g^n(\alpha(r)) = \alpha(r + nt)$ for any $r \in \mathbb{R}$. Thus by Theorem 2.2 (3) $\tau(g^n) = tn$.

The following proposition establishes a connection between translation numbers with respect to a word metric on a group and translation lengths of isometries of $CAT(0)$-spaces.

**Proposition 2.4.** *Let $G$ be a group and $S$ be a finite set of generators for $G$ closed under taking inverses. Suppose $G$ acts by isometries on a $CAT(0)$-space $(X, d)$ cocompactly ,i.e. the factor space $X/G$ is compact, and properly discontinuously, i.e. for any compact subset $K$ of $X$ the set $\{g \in G | g(K) \cap K \neq \emptyset\}$ is finite.*

*Then there is a constant $C > 0$ such that for any $g \in G$*

$$\tau(g)/C \leq \tau_S(g) \leq C\tau(g).$$

*Proof.* Pick a point $x_0 \in X$. We define a map $f: G \to X$ as follows:

$$f(g) = g(x_0)$$

for any $g \in G$.

By the Theorem 3.3.6 of [ECHLPT] there is a constant $C > 0$ such that for any $g, h \in H$

$$d(f(g), f(h))/C - C \leq d_S(g, h) \leq Cd(f(g), f(h)) + C$$

and for any $x \in X$ there is $g \in G$ such that $d(x, f(g)) \leq C$. If $g \in G$ is of finite order then $\tau(g) = \tau_S(g) = 0$. Suppose now that $g \in G$ has infinite order. Then by Theorem 2.2 (2) there is a biinfinite geodesic $\alpha: \mathbb{R} \to X$ such that for any $r \in \mathbb{R}$ $g(\alpha(r)) = \alpha(r + \tau(g))$, i.e. $g$ restricted to $\alpha$ is a translation by $\tau(g)$.

Put $x = \alpha(0)$.

Then by the triangle inequality

$$d(x_0, f(g^n)) = d(f(1), f(g^n)) \leq d(x_0, x) + d(x, g^n(x)) + d(g^n(x), g^n(x_0)) = 2d(x_0, x) + n\tau(g).$$

Analogously,

$$d(x_0, f(g^n)) \geq d(x, g^n(x)) - d(x_0, x) - d(g^n(x), g^n(x_0)) = -2d(x_0, x) + n\tau(g).$$

Thus

$$d_S(g^n, 1) \leq Cd(x_0, f(g^n)) + C \leq Cn\tau(g) + C(2d(x_0, x) + 1)$$

and

$$d_S(g^n, 1) \geq d(x_0, f(g^n))/C - C \geq n\tau(g)/C - 2d(x_0, x)/C - C$$

and therefore

$$\tau(g)/C \leq \tau_S(g) \leq C\tau(g)$$

and the proposition is proved.

Proposition 2.4 has the following easy but important corollary.

**Corollary 2.5.** *Let $G$, $S$ and $X$ be as in proposition 2.4. Then $G$ is strongly translation discrete.*

*Proof.* Recall that if $g$ is conjugate to $h$ in $G$ then $\tau(g) = \tau(h)$ (see Lemma 2.1(2)). If $g \in G$ is of infinite order then $\tau(g) > 0$ (see Theorem 2.2) and by Proposition 2.4 $\tau_S(g) > 0$. Thus $G$ is translation separable. Proposition 2.4 also implies that it suffices to show that for any $R > 0$ the set

$$\{[g] | \tau(g) \leq R\}$$

is finite (here $[g]$ is a conjugacy class of $g$).

Pick $r \geq 0$. Suppose $g \in G$ and $\tau(g) \leq R$. Since $X/G$ is compact, there is a compact subset $K$ of $X$ such that $GK = X$. Take any closed ball $B(y, R_0)$ with center $y$ and radius $R_0$ containing $K$. Then $B$ is compact since $X$ is proper and $GB = X$. Since $g$ is of infinite order, by Theorem 2.2 there is a geodesic $\alpha: \mathbb{R} \to X$ such that for any $r \in \mathbb{R}$ $g(\alpha(r)) = \alpha(r + \tau(g))$, i.e. $g$ restricted to $\alpha$ is a translation by $\tau(g)$.

Put $\alpha(0) = x_1$ and $t = \tau(g)$. Since $GB = X$, there is a point $x_0 \in B$ and $h \in G$ such that $h(x_0) = x_1$. Consider now a map $\beta: \mathbb{R} \to X$ defined as $\beta = h^{-1} \circ \alpha$. Then clearly $\beta$ is a geodesic and it has the property that for any $r \in \mathbb{R}$ $h^{-1}gh(\beta(r)) = \beta(r + t)$ and $\beta(0) = x_0$. Thus $d(x_0, h^{-1}gh(x_0)) = t \leq R$.

Clearly the set $B_1 = B(y, R_0 + R)$ is compact and $x_0, h^{-1}gh(x_0) \in B_1$. Thus we have shown that for any $g \in G$ with $0 < \tau(g) \le R$ there is an element $g_1$ conjugate to $g$ such that $g_1(B_1) \cap B_1 \ne \emptyset$. However the action of $G$ on $X$ is properly discontinuous and so the set

$$\{f \in G | f(B_1) \cap B_1 \ne \emptyset\}$$

is finite and therefore the set

$$\{[g] | \tau(g) \le R\}$$

is finite. This completes the proof of Corollary 2.5

We are now able to prove Theorem 0.1.

**Corollary 2.6.** *Let $G = <S|R>$ be a group presentation satisfying one of the C(4)-T(4)-P or C(3)-T(6)-P or C(6)-P conditions. Then $G$ is strongly translation discrete*

*Proof.* Suppose $G$ is as above. Let $K$ be the standard finite 2-complex associated to the given presentation of $G$ (see section 1 for definitions). Let $p \colon \tilde{K} \to K$ be the universal covering of $K$. It is not hard to see that, since the presentation for $G$ satisfies the $P$ small cancellation condition (every piece has length one), $\tilde{K}$ has a natural structure of a polygonal complex (in the sense of [BB]). Namely, put every 2-cell with $n$ 1-cells in its boundary in $\tilde{K}$ to be isometric to a regular euclidean $n$-gon of side 1 and every 1-cell to be isometric to a unit interval. Clearly the intersection of two distinct 2-cells in $\tilde{K}$ is either empty, or consists of a single vertex or is a single closed 1-cell.

Then to any path $p$ in $\tilde{K}$ which can be represented as a concatenation of a finite number of subpaths each of which lies inside a closed cell we can assign a length $L(p)$ as a sum of the lengths of these subpaths. We shall call such a $p$ a piecewise-cell path. We then define the distance $d(x,y)$ between points $x$ and $y$ of $\tilde{K}$ as the infimum of the lengths of all piecewise-cell paths connecting $x$ to $y$. Then (see Theorem 1.1 of [Br]) $(\tilde{K}, d)$ is a *complete geodesic metric space*. A detailed discussion about the properties of the metric space $(\tilde{K}, d)$ can be found in [Br].

There is a natural action of $G$ on $\tilde{K}$ by covering translations such that $\tilde{K}/G = K$. It is not hard to see that this action is properly discontinuous and cocompact. Moreover, any covering translation is an isometry of $\tilde{K}$ since this translation and its inverse preserve the lengths of piecewise-cell paths.

For any vertex $p$ of $\tilde{K}$ we define the link $D_p$ of $p$ as follows. $D_p$ is a graph whose set of vertices is the set of 2-cells of $\tilde{K}$ meeting $p$. Two vertices in $D_p$ are joined by an edge if the corresponding 2-cells in $\tilde{K}$ intersect by a one-cell. If these 2-cells in $\tilde{K}$ are an $n$-gon and an $m$-gon then we put the edge in $D_p$ to be isometric to an interval of length $\pi(n-1)/n + \pi(m-1)/m$. This defines a metric on each connected component of the graph $D_p$.

The small cancellation condition for $G$ ensures that for any vertex $p$ of $\tilde{K}$ the length of any nontrivial circuit in the link $D_p$ of $p$ has length at least $2\pi$. Therefore, by Theorem 15 of [B] $\tilde{K}$ is a non-positively curved space. Besides $\tilde{K}$ is simply connected and so by Theorem 14 of [B] any two points in $\tilde{K}$ can be joined by a unique geodesic segment. Therefore, by Theorem 7 of [B] $(\tilde{K}, d)$ is a $CAT(0)$-space. Thus by Corollary 2.6 the group $G$ is strongly translation discrete.

## 3. C(4)-T(4)-P-GROUPS

Let $G = <S|R>$ be a group presentation satisfying the $C'''(4) - T(4)$ condition.

It easily follows from [GS2, corollary of Theorem 2.1 for $(4,4)''$ disc diagrams] that for any letters $a, b \in X = S \cup S^{-1}$.

if $\overline{a} = \overline{b}$ then $a = b$ lexicographically,

if $\overline{a} = \overline{b^{-1}}$ then $a = b^{-1}$ lexicographically.

Moreover, if $w$ is a freely reduced word of length 1 or 2 in the alphabet $X$ then it is geodesic. We will rely on these facts throughout this section without further explanations.

The following definition is due to S.Gersten and H.Short [GS2].

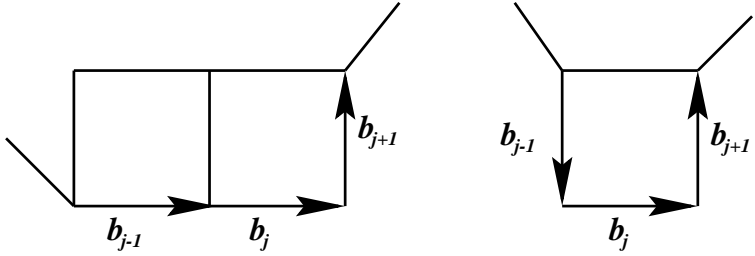

FIGURE 1

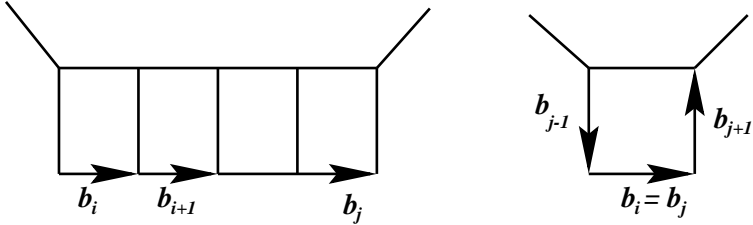

FIGURE 2

**Definition.** A word $w = b_0 b_1 \ldots b_n$ over the alphabet $X$ is called a $C''(4)$-geodesic if for all singular disc diagrams which includes $w$ as a part of their boundary cycles and for all $j$, $0 < j < n$, the pair of turning angles associated to $b_j$ does not consist of a non-negative angle followed by a positive angle (i.e. pairs $(0, \pi/2)$ and $(\pi/2, \pi/2)$ are forbidden) (See fig.1).

The following proposition is also proved in [GS2].

**Proposition 3.1.** *Let $G = <S|R>$ be a group presentation satisfying the $C''(4) - T(4)$ condition and let $X = S \cup S^{-1}$ define the word metric $d_X$ on $G$. Then*
  (a) *any $C''(4)$-geodesic word is $d_X$-geodesic;*
  (b) *any element $g \in G$ has a $C''(4)$-geodesic representative;*
  (c) *the set of all $C''(4)$-geodesic words is regular and form an biautomatic structure for $G$.*

We need also another notion from [ASc] and [GS2]. A subword $b_i \ldots b_j$, $0 < i \leq j < n$, of a word $w = b_0 b_1 \ldots b_n$ is called a **chain** if there is a disc diagram $D$ which includes $w$ as a part of its boundary cycle such that the first turning angle associated to $b_i$ and the second turning angle associated to $b_j$ are $\pi/2$, and all other turning angles between them are zero (see fig.2). If $i = j$ then the chain is said to be **degenerate**. An edge-path on the boundary of $D$ which corresponds to the chain $b_i \ldots b_j$ with its first and last vertices removed is termed an **interior** of this chain.

The following lemma gives a description of all geodesic words in the alphabet $X$. It is not hard to deduce it from Lemma 2 proved by K.Appel and P.Schupp in [ASc] or Corollary 5.1 in [GS2]. We shall present a proof of this fact for the sake of completeness.

**Lemma 3.2.** *Let $w = x_1 \ldots x_n$ be a word in the alphabet $X$. Then $w$ is geodesic if and only if it is freely reduced and does not contain a bad subword, i.e. a subword $x_i x_{i+1} \ldots x_j$ which can be included in a diagram as in fig.3*

*Proof.* Certainly if $w$ is not freely reduced or it contains a bad subword then it is not geodesic. So we need to prove the statement of the lemma only in the opposite direction.

We will proceed by induction on $l(w)$ - the length of $w$.

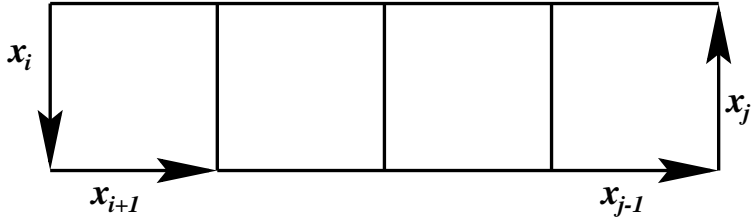

FIGURE 3

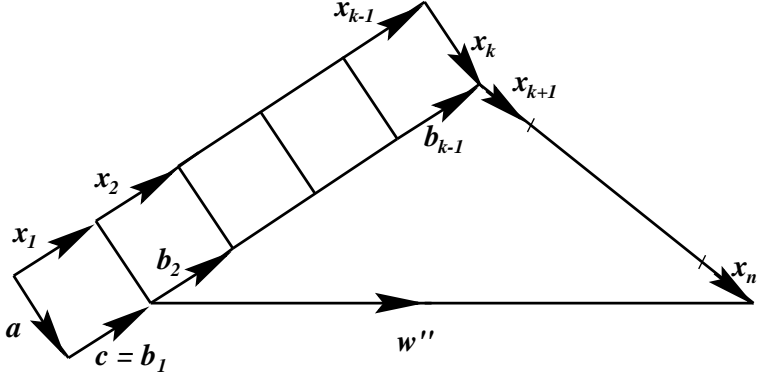

FIGURE 4

For $l(w) \leq 3$ the statement easily follows from [GS2, Corollary of Theorem 2.1 for $(4,4)''$ diagrams].

Suppose we have already proved Lemma 3.2 for words of length less than $n$. Suppose there is a freely reduced word $w$ of length $n$ which does not contain bad subwords and which is not geodesic.

Let $w'$ be a $C''(4)$ geodesic word such that $\overline{w} = \overline{w'}$. In particular $l(w') < l(w)$.

There is a singular disc diagram $D$ for the word $w'w^{-1}$ representing the identity element in $G$.

If this diagram has a cut vertex then $w = w_1w_2$, $w' = w'_1w'_2$ such that $\overline{w_i} = \overline{w'_i}$ for $i = 1,2$.

Since $l(w') < l(w)$ we have $l(w'_i) < l(w_i)$ for some $i$. Since $w_i$ is obviously freely reduced, not geodesic and has length less than $n$ we may apply the inductive hypothesis to $w_i$ and conclude that it contains a bad subword. But this bad subword is also a subword of $w$, which contradicts our assumptions.

So we can assume that our diagram does not have a cut vertex, i.e. this is a PE-disc diagram.

Then by [GS2, Corollary 5.1] (or [ASc]) on $\partial D$ there are four chains with disjoint interiors. Since $w'$ is geodesic and $w$ does not contain bad subwords, any of these chains cannot sit strictly inside $w$ or $w'$. So the only possibility is that two of these chains have an initial vertex $v_0$ of $w'$ in common and the other two have the last vertex $v_1$ of $w'$ in common. Moreover, the chain $c_1$ which meets $v_0$ and whose interior lies iside $w'$ must be degenerate because otherwise the word $w'$ fails to be $C''(4)$-geodesic.

So we have a picture like the one in fig.4 where $w = x_1x_2\ldots x_n$ and $w' = ab_1w''$.

The word $w''$ is a subword of $w'$ and so it is geodesic and has length $l(w'') = l(w') - 2$. Also we know that the word

$$w_1 = b_2b_3\ldots b_{k-1}x_{k+1}\ldots x_n$$

represents the same element as $w''$ and has length $n - 2 = l(w) - 2$. It is easy to see that $w_1$ cannot be geodesic because otherwise $l(w') - 2 = l(w) - 2$ and $w$ fails to be not geodesic.

Besides the word $b_2b_3\ldots b_{k-1}$ is freely reduced because the word $x_2x_3\ldots x_{k-1}$ is freely reduced and the $T(4)$ condition holds. So in order to see that $w_1$ is freely reduced it is sufficient to check that $b_{k-1} \neq x_{k+1}^{-1}$.

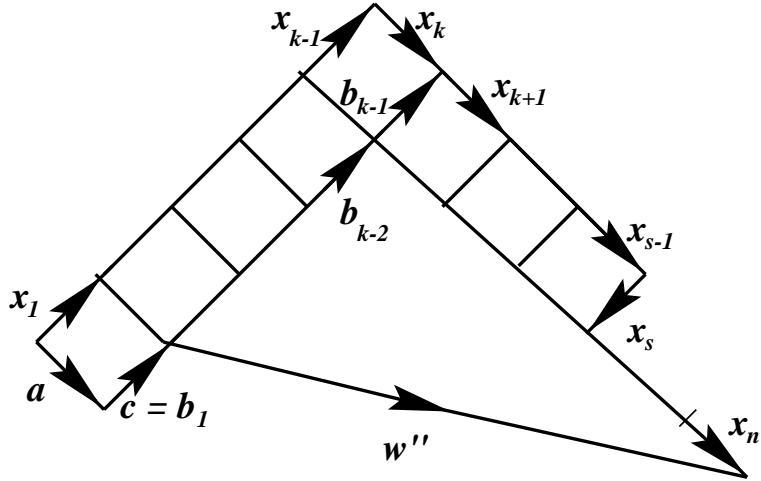

FIGURE 5

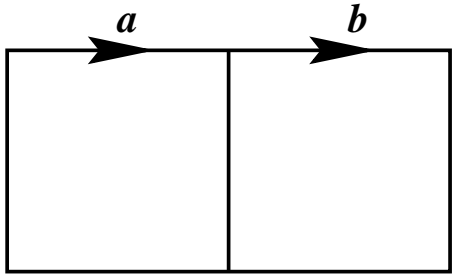

FIGURE 6

Suppose $b_{k-1} = x_{k+1}^{-1}$.

Then the word $w$ itself contains a forbidden subword $x_{k-1}x_k x_{k+1} = x_{k-1}x_k b_{k-1}^{-1}$ which is a 3/4-relator what contradicts our assumptions about $w$.

Thus the word $w_1$ is freely reduced, not geodesic and by the inductive hypothesis it contains a bad subword. This bad subword cannot lie wholly inside the word $x_{k+1}\ldots x_n$. Besides the $T(4)$ condition tells us that $b_i b_{i+1}$ cannot be a half-relator and therefore our bad subword cannot contain more than one letter of $b_2\ldots b_{k-1}$. So the only possibility is that this bad subword has the form $b_{k-1}x_{k+1}\ldots x_{s-1}x_s$. But in this case the word $x_{k-1}x_k x_{k+1}\ldots x_{s-1}x_s$ is a bad word (see fig.5) and it is a subword of $w$.

This contradiction completes the proof of Lemma 2.2.

*Remark 3.3.* Note that for $C''(4) - T(4)$ presentations if for letters $a$ and $b$ such that $a \ne b^{-1}$ there is a diagram as in fig.6 then $ab$ is not a half-relator.

**Corollary 3.4.** *Let $G = < S|R >$ be a group presentation satisfying the $C''(4) - T(4)$ condition and $X = S \cup S^{-1}$. Then the set of all $d_X$-geodesic words is a regular language.*

*Proof.* This is an immediate corollary of Lemma 3.2 and the fact that if a freely reduced word $w = x_1\ldots x_n$ can be encluded in a diagram $D$ as in fig.7 then this diagram is unique (because of the T(4)-P-condition).

We shall give an informal description of the automaton $M$ recognizing the set of all geodesic words. It rejects all non freely reduced words. Reading a word it looks for pairs of consecutive letters which form a half-relator. When such a pair is found, $M$ begins to build (if possible) a diagram corresponding to a potentially bad subword beginning with that half-relator. Analysing a current letter, $M$ needs to know not the whole initial segment of an already built potentially bad diagram but only the last relator of it (see the remark in the beginning of the proof). If $M$ succeds in building a bad subword diagram then it rejects the word without looking at the rest of it. If $M$ finds a half-relator which does not end the bad diagram then $M$ starts a new building process beginning from this half-relator. If at some moment $M$ is unable to continue building the bad diagram by a reason other than mentioned above, it gives up the current building process and starts a new search for a half-relator. If after reading the whole input word $M$ does not find a bad subword, it accepts this word. It is not hard to see that the language recognized by this automaton is the language of all geodesics.

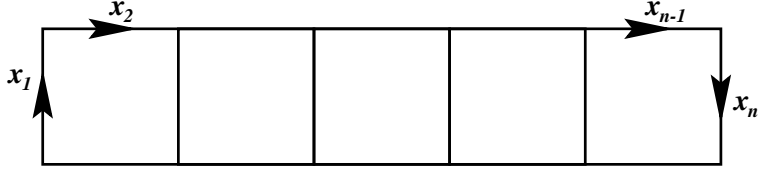

FIGURE 7

**Proposition 3.5.** Let $G = <S|R>$ be a $C''(4) - T(4)$ presentation and $X = S \cup S^{-1}$ a generating set defining the word metric $d_X$ on $G$. Let $\alpha$ be a conjugacy class of length at least three and let a $d_X$-geodesic word $w$ represent some shortest element of $\alpha$.

Then either $w^k$ is geodesic for any integer $k \geq 1$ or there is a cyclic permutation $u$ of $w$ such that

$$d_X(\overline{u^k}, 1) = k \cdot (l(w) - 1)$$

for even $k \geq 1$ and

$$d_X(\overline{u^k}, 1) = k \cdot (l(w) - 1) + 1$$

for odd $k \geq 1$. In this last case $g = \overline{w^2}$ is periodically geodesic, that is $l_X(g^n) = |n|l_X(g)$ for each $n \in \mathbb{Z}$.

*Proof.*

Let $w$ be as in Proposition 3.5. If $w^k$ is always geodesic we have nothing to prove.

Suppose there is an integer $k > 1$ such that $w^k$ is not geodesic. Since $w$ represents a shortest element in the conjugacy class $\alpha$ of $\overline{w}$ and $w$ is geodesic then any cyclic permutation of $w$ is also geodesic. Thus $w^k$ is locally $n$-geodesic $n = l(w)$ and in particular freely reduced. (We say that a word $y$ in the alphabet $X$ is $m$-locally geodesic if any subword $z$ of $y$ of length at most $m$ is $d_X$-geodesic.)

So by Lemma 3.2 the word $w^k$ contains a bad subword $v$ of length $s \geq 3$. Since $w^k$ is locally $n$-geodesic, we have $s > n$. On the other hand, $s \leq n + 2$. To see this suppose our bad word has length $s \geq n + 3$. Then any two-letter subword $ab$ of $w^k$ may be included in a diagram as in fig.6 and therefore by Remark 3.3 cannot be a half-relator. However, the first two letters of our bad word $v$ give a half-relator what is impossible.

Thus we have established that $n + 1 \leq s \leq n + 2$. This means that for some cyclic permutation $u = a_1 a_2 \ldots a_n$ of $w$ one of the following two cases holds.

*Case 1.*

There is a diagram as one on the fig.8.

*Remark 3.6.* For $C''(4) - T(4)$ presentations if $a \neq b^{-1}$ and $ab$ is a half of some relator then this relator is unique.

So by Remark 3.6 the first and the last squares of the diagram in fig.8 correspond to the same relator and we indeed have a picture labeled like one in fig.8. Therefore

$$\overline{a_2 a_3 \ldots a_n a_1} = \overline{c^{-1} b_2 b_3 \ldots b_{n-1} c}$$

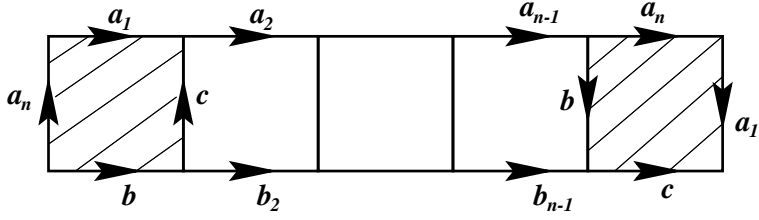

FIGURE 8

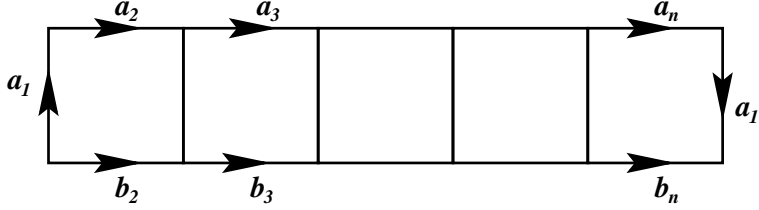

FIGURE 9

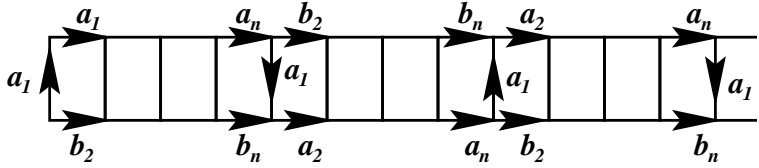

FIGURE 10

and we have found in the conjugacy class $\alpha$ an element of length less than $n$. This contradicts our assumptions about minimality of $w$.

*Case 2.*

A diagram like one in figure 9 exists.

It is easy to see that the word $b_2 b_3 \ldots b_n$ is freely reduced since the word $a_2 a_3 \ldots a_n$ is freely reduced.

We will use this diagram as a "brick" to build a long wall as shown in fig.10

Suppose $a_n = b_2^{-1}$. Then $b_2^{-1} a_1 a_2 \equiv a_n a_1 a_2$ is a 3/4-relator and the word $w^k$ is not locally 3-geodesic what is impossible since $l(w) \geq 3$ by our assumptions.

Analagously, $b_n \neq a_2^{-1}$ since otherwise $a_n a_1 b_n^{-1} \equiv a_n a_1 a_2$ is a 3/4-relator.

Thus we see that lower and upper borders of our wall represent freely reduced and, moreover, geodesic (by Lemma 3.2) words. Clearly

$$\overline{(a_1 a_2 \ldots a_n)^{2k}} = \overline{(b_2 b_3 \ldots b_n a_2 a_2 \ldots a_n)^k}$$

and the latter word is geodesic.

Analogously,

$$\overline{(a_1 a_2 \ldots a_n)^{2k+1}} = \overline{a_1 a_2 \ldots a_n (b_2 b_3 \ldots b_n a_2 a_2 \ldots a_n)^k}$$

and the latter word is geodesic by Lemma 3.2

Thus we have found that for cyclic permutation $u = a_1 a_2 \ldots a_n$ of $w$, for odd $k$

$$d_X(\overline{u^k}, 1) = 1 + k(n-1),$$

and for even $k$
$$d_X(\overline{u^k}, 1) = k(n-1).$$

It is also clear by Lemma 3.2 that the word $v = b_2 \ldots b_n a_2 \ldots a_n$, representing $\overline{w^2}$, is periodically geodesic. This completes the proof of the proposition.

**Corollary 3.7.** *If $G = <S|R>$ is a group presentation satisfying the $C''(4) - T(4)$ condition and $X = S \cup S^{-1}$ then for any congugacy class $\alpha$ of length $d_X(\alpha) \geq 3$ we have $\tau(\alpha) = l_X(\alpha)$ or $l_X(\alpha) - 1$. This means that the group $G$ is strongly translation descrete (It is translation separable because it is biautomatic.)*

**Corollary 3.8.** *If $G = <S|R>$ is a group presentation satisfying the $C(4) - T(4) - P$ condition then $G$ is strongly translation discrete.*

*Proof.*

Indeed, as Gersten and Short pointed out in [GS2], there is a free group of finite rank $F$ such that the group $G_1 = G * F$ admits a $C''(4) - T(4)$ presentation and, therefore, is strongly translation discrete.

Choose generating sets $X, Y, Z = X \cup Y$ for groups $G, F, G_1$ accordingly in such a way that $d_X|_G = d_Z|_{G_1}$. Denote for any element $g \in G$ its conjugacy class in $G$ by $[g]_G$ and its conjugacy class in $G_1$ by $[g]_{G_1}$. Clearly $[g]_G = [g]_{G_1} \cap G$ for any $g \in G$. So for any real $r$

$$|\{[g]_G | \tau_X([g]_G) \leq r, g \in G\}| \leq$$
$$\leq |\{[g]_{G_1} | \tau_Z([g]_{G_1}) \leq r, g \in G\}| \leq$$
$$\leq |\{[g_1]_{G_1} | \tau_Z([g_1]_{G_1}) \leq r, g_1 \in G_1\}| < \infty$$

This completes the proof of Corollary 3.8.

*Remark 3.9.* An argument analogous to the one used in the proof of Proposition 3.5 shows that if $G = <S|R>$ is a group presentation satisfying the $C''(4) - T(4)$ condition and $\alpha$ is a conjugacy class, then

$$l_X(\alpha) = 1 \Rightarrow \tau_X(\alpha) = 1$$

and

$$l_X(\alpha) = 2 \Rightarrow \tau_X(\alpha) = 1 \text{ or } \tau_X(\alpha) = 2.$$

Also, if $l_X(\alpha) \leq 2$ and $g \in \alpha$ is a shortest element then $g^2$ is periodically geodesic.

**Corollary 3.10.** *Let $G$ be a group as in Corollary 3.8. Then for any $g \in G$ of infinite order $\sup\{n | g = x^n \text{ for some } x \in G\} < \infty$.*

In particular $G$ does not contain subgroups isomorphic to the group of $p$-adic fractions $\mathbb{Q}_p = \{k/p^l | k \in \mathbb{Z}, l \in \mathbb{N}\}$ (where $p$ is some prime number) and therefore to the group of rational numbers $\mathbb{Q}$.

*Proof.* It immediately follows from property (d) of translation numbers given in the introduction and the fact that $G$ is translation discrete and so 0 is not an accumulation point for the set of translation numbers.

**Corollary 3.11.** *If $G = <S|R>$ is a group presentation satisfying the $C(4) - T(4) - P$ condition and $X = S \cup S^{-1}$ then for any freely reduced word $w$ in the alphabet $X$ there is an algorithm which determines whether or not the equation*

$$x^n = \overline{w} \tag{1}$$

*has a solution in $G$. Moreover, for any word $w$ there is an algorithm for calculating $\max\{n \mid x^n = \overline{w} \text{ for some } x \in G\}$.*

*Proof.* It follows from Remark 3.9 and Corollary 3.7 that $C(4)'' - T(4)$ and therefore $C(4) - T(4) - P$ groups are torsion free. As proved in [GS2, Proposition 5.9], $G$ can be effectively embedded in a group $G_1$ given by a $C''(4) - T(4)$ presentation and, moreover, the image of $G$ is a free factor in $G_1$. Thus it is not hard to see that for $g \in G$ the equation $x^n = g$ has a solution in $G$ if and only if it has a solution in $G_1$. Therefore we can assume from the beginning that our group $G$ is given by a $C(4)'' - T(4)$ presentation. As shown in [GS2] in this situation one can effectively construct a biatomatic structure with uniqueness

$$(L, M_=, R_{x_1}, .., R_{x_n}, L_{x_1}, .., L_{x_n})$$

on $G$ where $L$ is a sublanguage of the set of all geodesic words, $M_=$ is the equality checker and $R_{x_i}$, $L_{x_i}$ are automata recognizing right and left translations by the letter $x_i$. Therefore, as shown in [GS3, Theorem 8.3], one can use this biautomatic structure to find a word $u \in L$ such that $\overline{u}$ is a shortest element in the conjugacy class of $\overline{w}$. If $u$ is the empty word then $\overline{w} = 1$ and (1) has a solution. If $l(u) < n$ and $u$ is not empty then (1) has no solution. Indeed, if $h^n = \overline{w}$ for some $h$ then $f^n = \overline{u}$ for some $f \neq 1$. Thus $\tau_X(f) = \tau_X(u)/n < 1$ which is impossible as Corollary 3.7 and Remark 3.9 show. If $u$ is not empty and $l(u) \geq n$ then for any $f$ such that $f^n = \overline{u}$ we have

$$l(u)/n \geq \tau_X(u)/n = \tau_X(f) \geq l_X(\alpha) - 1$$

where $\alpha$ is the congugacy class of $f$. Therefore now one can check (using the biautomatic structure on $G$) if there are words $v \in L$ such that $l(v) \leq (l(u)/n) + 1$ and the element $\overline{v^n}$ is conjugate to $\overline{u}$. If there is such $v$ then (1) has a solution and and if there is no such $v$ then (1) has no solution. Thus for any $n$ we have constructed an algorithm which determines wheter or not (1) has a solution in $G$.

As we have seen if $l(u) < n$ and $u$ is not empty then (1) has no solution in $G$. So to calculate $\max\{n | x^n = \overline{w}$ for some $x \in G\}$ it suffices to check those $n$ for which $n \leq l(u)$. This yields the second part of Corollary 3.11.

## 4. C(3)-T(6)-P AND C(6)-P GROUPS.

The situation with $C(3) - T(6) - P$ and $C(6) - P$ groups is essentially the same as with $C(4) - T(4) - P$ groups, so we shall give just brief sketches of proofs.

First, we establish strong translation discreteness for $C'''(3) - T(6)$ groups. Since (see [GS2]) $C(3) - T(6) - P$ and $C(6) - P$ groups occur as free factors of $C(3) - T(6) - P$ and $C(6) - P$ groups, the general case follows as in Corollary 3.8.

Let $G = <S|R>$ be a group presentation satisfying the $C'''(3) - T(6)$ condition.

It easily follows from [GS2, corollary of Theorem 2.1 for $(3,6)''$ disc diagrams] that for any letters $a, b \in X = S \cup S^{-1}$

if $\overline{a} = \overline{b}$ then $a = b$ lexocagraphically,

if $\overline{a} = \overline{b^{-1}}$ then $a = b^{-1}$ lexocagraphically.

Moreover, all letters of the alphabet $X$ are $d_X$-geodesic words. We will rely on these facts throughout this section without further explanations.

**Definition.** A word $w = b_0 b_1 \ldots b_n$ over the alphabet $X$ is called a $C'''(3)$-geodesic if for all singular disc diagrams which include $w$ as a part of their boundary cycles and for all $j$, $0 < j < n$, the pair of turning angles $(\theta, \theta')$ associated to $b_j$ does not consist of a non-negative angle followed by a positive angle and both $\theta$ and $\theta'$ are different from $2\pi/3$.

The following proposition is also proved in [GS2].

**Proposition 4.1.** *Let $G = <S|R>$ be a group presentation satisfying the $C'''(3) - T(6)$ condition and $X = S \cup S^{-1}$ defines the word metric $d_X$ on $G$. Then*
 (a) *any $C'''(3)$-geodesic word is $d_X$-geodesic;*
 (b) *any element $g \in G$ has a $C'''(3)$-geodesic representative;*
 (c) *the set of all $C'''(3)$-geodesic words is regular and form an biautomatic structure for $G$.*

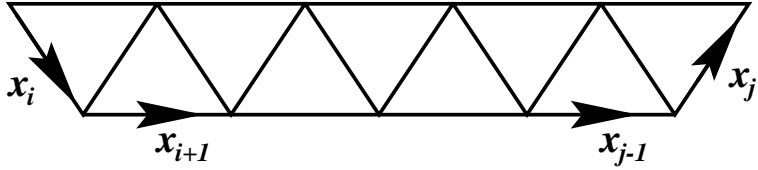

FIGURE 11

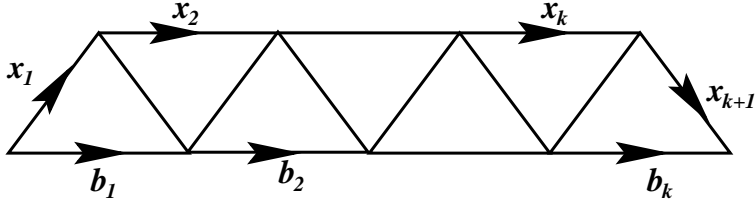

FIGURE 12

We need also another notion from [GS2]. A subword $b_i \ldots b_j$, $0 < i \leq j < n$, of a word $w = b_0 b_1 \ldots b_n$ is called a **chain** if there is a disc diagram $D$ which includes $w$ as a part of its boundary cycle such that the first turning angle associated to $b_i$ and the second turning angle associated to $b_j$ are positive, and all other turning angles between them are zero. An edge-path in the boundary of $D$ which corresponds to the chain $b_i \ldots b_j$ with its first and last vertices removed is termed an **interior** of this chain. Note that a chain may consist of a single letter.

As in the previous section at first we give a description of all geodesic words in the alphabet $X$.

**Lemma 4.2.** Let $w = x_1 \ldots x_n$ be a word in the alphabet $X$. Then $w$ is geodesic if and only if it is freely reduced and does not contain a bad subword, i.e. a subword $x_i x_{i+1} \ldots x_j$ which can be included in a diagram as in fig.11.

*Proof.* The proof proceeds by induction on $l(w)$. As before we suppose that there is a freely reduced nongeodesic word $w$ without bad subwords and consider a $C''(3)$ geodesic word $\underline{w'}$ representing the same element as $w$. Consider a singular disc diagram $D$ corresponding to the relation $\overline{ww'^{-1}} = 1$. If this diagram has a cut vertex then the statement follows by induction. If there are no cut vertices then $D$ is a PE-disc diagram and by [GS2, Corollary 6.1] on the boundary of $D$ one of the following occurs

1) there are three vertices with turning angle $2\pi/3$,
2) there are two vertices with turning angle $2\pi/3$ and four chains with disjoint interiors,
3) there are one vertex with turning angle $2\pi/3$ and five chains with disjoint interiors,
4) there are six chains with disjoint interiors.
It is easy to see that any of these cases gives us a contradiction.

The following observations immediately follow from the previous lemma and the $C''(3) - T(6)$ condition.

*Remark 4.3..* If a freely reduced word $x_1 x_2 \ldots x_k x_{k+1}$ can be included in the diagram as in fig.12 then $b_1 b_2 \ldots b_k$ is freely reduced and geodesic and words $x_1 x_2 \ldots x_k$ and $x_2 x_3 \ldots x_k x_{k+1}$ are also geodesic.

*Remark 4.4..* If for two letters $a$ and $b$ the word $ab$ is a 2/3 of some relator, then this relator is unique. If for two non mutually inverse letters $x_1$ and $x_2$ there is a diagram $D$ as in fig.13 then this diagram is unique.

The following statement is analogous to Proposition 2.5.

**Proposition 4.5.** Let $G = <S|R>$ be a $C''(3) - T(6)$ presentation and let $X = S \cup S^{-1}$ be a generating set defining the word metric $d_X$ on $G$. Let $\alpha$ be a conjugacy class of length at least two and let a geodesic word $w$ represent some shortest element of $\alpha$.

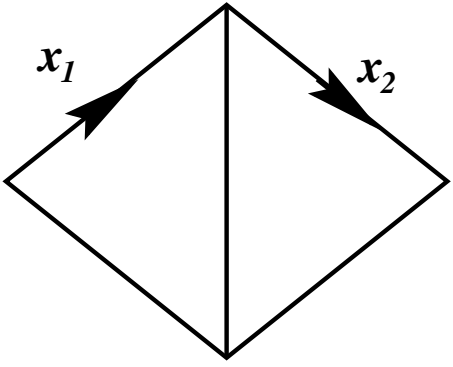

FIGURE 13

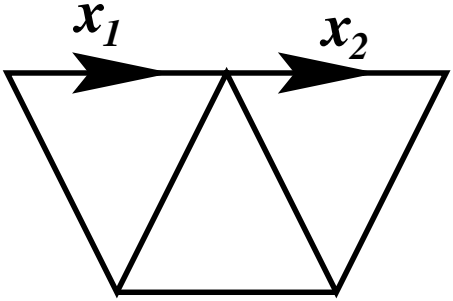

FIGURE 14

Then either $w^k$ is geodesic for any integer $k \geq 1$ or there is a cyclic permutation $u$ of $w$ such that

$$d_X(\overline{u^{2k}}, 1) = 2k \cdot l(w) - k$$

and

$$d_X(\overline{u^{2k+1}}, 1) = (2k+1) \cdot l(w) - k$$

for any integer $k \geq 0$. Also, $g = \overline{w^2}$ is periodically geodesic that is $l_X(g^n) = |n| l_X(g)$ for each $n \in \mathbb{Z}$.

*Proof.*

As before any cyclic permutation of $w$ is geodesic because of minimality of $d_X(\overline{w}, 1)$ in the conjugacy class $\alpha$ of $\overline{w}$ and therefore the word $w^k$ is locally $l(w)$-geodesic for any $k \geq 1$. If $w^k$ is geodesic for each $k$, the proposition holds.

If $w^k$ is not geodesic for some $k$ then by Lemma 4.2 it contains a bad subword of length $s$. Clearly $s \geq l(w) + 1$ since the word $w^k$ is locally $l(w)$-geodesic. On the other hand $s \leq l(w) + 2$. Otherwise any two-letter subword $x_1 x_2$ of $w$ and any its cyclic permutation can be included in a diagram like the one in fig.14 . Therefore by the T(6)-condition $x_1 x_2$ cannot be included in a diagram as in fig.13 which contradicts the fact that such two-letter words should occur in the beginning of our bad subword.

Thus for some cuclic permutation $u = a_1 \ldots a_l$ of $w$ one of the following holds.

*Case* 1. There is a diagram like the one in fig.15.

Then by Remark 4.4 the first and the last paralellograms in fig.15 must correspond to the same diagram. So we have the equality

$$\overline{a_1 a_2 \ldots a_l} = \overline{d b_2 \ldots b_l d^{-1}}$$

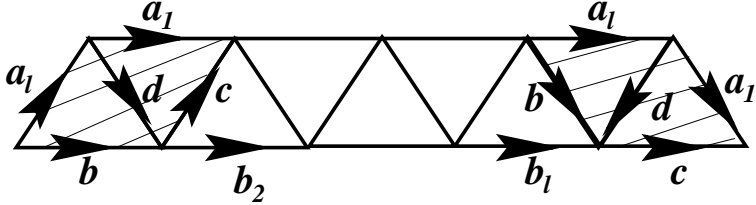

FIGURE 15

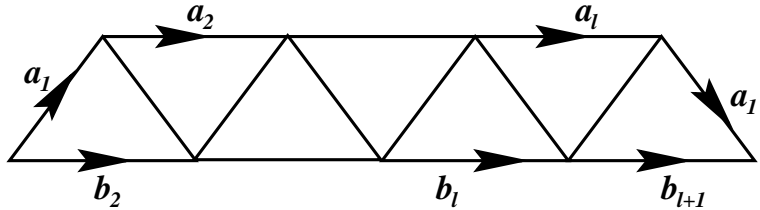

FIGURE 16

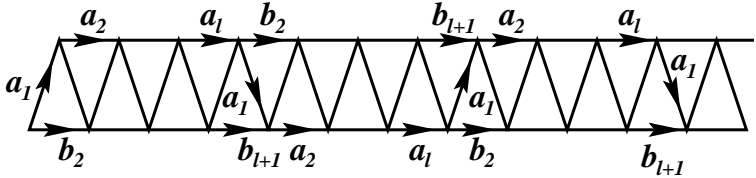

FIGURE 17

which contradicts the minimality of $\overline{w}$ in its conjugacy class.

$Case\,2.$ There is a diagram like the one in fig.16.

Using this diagram as a "brick", we can build an arbitrary long "wall" (see fig.17).

It is not hard to see that upper and lower boundaries of it corespond to freely reduced words. Moreovere, it follows from Lemma 4.2 that the words

$$a_1(a_2\ldots a_l b_2\ldots b_{l+1})^k$$

and

$$(b_2\ldots b_{l+1} a_2\ldots a_l)^k$$

are geodesic. It is also clear from Lemma 4.2 that the word $v = b_2\ldots b_{l+1} a_2 \ldots a_l$, representing $\overline{w^2}$, is periodically geodesic which implies the statement of the proposition.

**Corollary 4.6.** *If $G = <S|R>$ is a group presentation satisfying the $C''(3) - T(6)$ condition and $X = S \cup S^{-1}$ then for any congugacy class $\alpha$ of length $l_X(\alpha) \geq 2$ we have $\tau(\alpha) = l_X(\alpha)$ or $l_X(\alpha) - 1/2$. This means that the group $G$ is strongly translation discrete.*

**Corollary 4.7.** *If $G = <S|R>$ is a group presentation satisfying the $C(3) - T(6) - P$ or $C(6) - P$ condition then $G$ is strongly translation descrete.*

*Proof.* As we pointed out in the beginning of this section, the statement follows exactly as in Corollary 3.8 from Corollary 4.6 and the fact that all $C(3) - T(6) - P$-groups and $C(6) - P$-groups are free factors of $C''(3) - T(6)$ groups.

*Remark 4.8.* It can be established by an argument analogous to the one used in the proof of Proposition 4.5 that if $\alpha$ is a conjugacy class of length one then either $\tau_X(\alpha) = 1$ or $\tau_X(\alpha) = 1/2$. Moreover, if $g \in \alpha$ is a shortest element then $g^2$ is periodically geodesic.

**Corollary 4.9.** *If If $G = <S|R>$ is a group presentation satisfying the $C(3) - T(6) - P$ or $C(6) - P$ condition and $X = S \cup S^{-1}$ then for any freely reduced word $w$ in the alphabet $X$ there is an algorithm which determines whether or not the equation*

$$x^n = \overline{w} \qquad (2)$$

*has a solution in $G$. Also, one can determine algorithmically the maximum $n$ such that $w$ is an $n$-th power in $G$.*

*Proof.* It follows from Remark 4.8 and Corollary 4.6 that $C(3)'' - T(6)$ and therefore $C(3) - T(6) - P$ groups are torsion free. As proved in [GS2, Section 6], $G$ can be effectively embedded in a group $G_1$ given by a $C''(3) - T(6)$ presentation and, moreover, the image of $G$ is a free factor in $G_1$. Thus it is not hard to see that for $g \in G$ equation (2) has a solution in $G$ if and only if it has a solution in $G_1$. Therefore we can assume from the beginning that our group $G$ is given by a $C(3)'' - T(6)$ presentation. As shown in [GS2] in this situation one can effectively construct a biatomatic structure with uniqueness

$$(L, M_=, R_{x_1}, .., R_{x_n}, L_{x_1}, .., L_{x_n})$$

on $G$ where $L$ is a sublanguage of the set of all geodesic words, $M_=$ is the equality checker and $R_{x_i}, L_{x_i}$ are automata recognizing right and left translations by the letter $x_i$. Therefore, as shown in [GS3, Theorem 8.3], one can use this biautomatic structure to find a word $u \in L$ such that $\overline{u}$ is a shortest element in the conjugacy class of $\overline{w}$. If $u$ is an empty then $\overline{w} = 1$ and (1) has a solution. If $l(u)/n < 1/2$ and $u$ is not empty then (1) has no solution. Indeed, if $h^n = \overline{w}$ for some $h$ then $f^n = \overline{u}$ for some $f \neq 1$. Thus $\tau_X(f) = \tau_X(u)/n < 1/2$ which is impossible as Corollary 4.6 and Remark 4.8 show. If $u$ is not empty and $l(u)/n \geq 1/2$ then for any $f$ such that $f^n = \overline{u}$ we have

$$l(u)/n \geq \tau_X(u)/n = \tau_X(f) \geq l_X(\alpha) - 1/2$$

where $\alpha$ is the congugacy class of $f$. Therefore now one can check (using the biautomatic structure on $G$) if there are words $v \in L$ such that $l(v) \leq (l(u)/n) + 1/2$ and the element $\overline{v^n}$ is conjugate to $\overline{u}$. If there is such $v$ then (2) has a solution and and if there is no such $v$ then (2) has no solution.

Thus for any $n$ we have constructed an algorithm which determines wheter or not (2) has a solution in $G$.

As we have seen if $l(u)/n < 1/2$ and $u$ is not empty then (2) has no solution in $G$. So to calculate $\max\{n | x^n = \overline{w} \text{ for some } x \in G\}$ it suffices to check those $n$ for which $n \leq 2l(u)$. This yields the second part of Corollary 4.9.

*Remark 4.10.* It is not hard to see that, proceeding as in the proof of Corollary 3.4, one can show that for a group $G$ given by a $C''(3) - T(6)$ presentation $G = <S|R>$ the language of all geodesic words in $X = S \cup S^{-1}$ is regular.

*Remark 4.11.* Corollary 4.9 and Corollary 3.11 yield Theorem 0.4 and Theorem 0.5 from the Introduction. It also follows from our Corollary 3.7, Remark 3.9, Corollary 4.6 and Remark 4.8 that if $G = <S|R>$ is a $C''(4) - T(4)$ or $C''(3) - T(6)$ presentation and $X = S \cup S^{-1}$ then for any $g \in G$ the number $\tau_X(g)$ is rational and, moreover, $2\tau_X(g)$ is integer. Also, the square of every nontrivial element is conjugate to a periodically geodesic element. This gives us the statement of Theorem 0.7 from the Introduction.

## 5. Open questions.

There are a lot of open questions connected with translation numbers. (In fact most questions fomulated by S.Gersten and H.Short in [GS1] are still not resolved and can be reformulated in the language of translation discreteness.) We shall point out here only a few of them.

**Question 1.** Are automatic groups translation separable? Do biautomatic (automatic) groups have some stronger kind of translation discreteness?

It seems that in full generality these questions are very difficult. However a lot of interesting classes of groups, such as geometrically finite groups, fundamental groups of Haken 3-manifolds, braid groups, some more general classes of small cancellation groups than those we have considered here, etc., were proved to be automatic (biautomatic). It would be interesting to find out what kind of translation discreteness these groups posess.

**Question 2.** Is it true that for automatic (biautomatic) group the number of conjugacy classes of elements of finite order is always finite? (This property would follow from being strongly translation discrete.)

**Question 3.** Can an automatic (biautomatic) group contain a subgroup isomorphic to the additive group of rational numbers (p-adic fractions)? (The negative answer would follow from being translation discrete.)

**Question 4.** Can a translation number in an automatic group be irrational? What about a group acting cocompactly and properly discontinuously on a CAT(0)-space? It is known (see, for example, [BGSS,proof of Theorem D]) that if $G$ is hyperbolic then $\tau_X(g)$ is rational for any $g \in G$ and for any finite generating set $X$ closed under taking inverses. Moreover, Theorem 13 of [S] implies that in this situation there is an integer $M \neq 0$ such that $M \cdot \tau_X(G) \subseteq \mathbb{Z}$. Translation numbers for finitely generated nilpotent groups are also rational [C2]. Our Theorem 0.7 shows that if $G = <S|R>$ is a group presentation satisfying $C''(4) - T(4)$ or $C''(3) - T(6)$ small cancellation condition and $X = S \cup S^{-1}$ then $2\tau_X(G) \subseteq \mathbb{Z}$. On the other hand, in [C3] G.Conner constructs a group which has elements with irrational translation numbers.

**Question 5.** Suppose $G$ has a C(4)-T(4)-P, C(6)-P or C(3)-T(6)-P presentation. Is the power-conjugacy problem solvable in $G$, that is, given two elements of $G$, is there a way to determine whether some powers of these elements are conjugate?

Corollary 0.6 seems to indicate that the answer is positive. Also, M.Anshel and P.Stebe [AS] and L.Comerford [Com2] showed that the power-conjugacy problem is solvable for some other classes of small cancellation groups and their cyclic amalgamations.

## References

[AS]    M.Anshel and P.Stebe, *Conjugate powers in free products with amalgamations*, Houston J. Math. **2** (1976), no. 2, 139–147.

[ASc]   K.Appel and P.Schupp, *Artin groups and infinite Coxeter groups*, Invent. Math. **72** (1983), 201-220.

[B]     W.Ballman, *Singular spaces of non-positive curvature*, Sur les Groupes Hyperboliques d'apres Mikhael Gromov (E.Ghys and P. de la Harpe, eds.), Birkhauser, Berlin, 1990 pages 189–201.

[BB]    W.Ballmann and M.Brin, *Polygonal comlexes and combinatorial group theory*, preprint, University of Maryland at College Park, 1992.

[BGrS]  W.Ballman, M.Gromov, V.Schroeder, *Manifolds of nonpositive curvature*, Progress in Mathematics, vol. 61, Birhauser, Boston, 1985.

[BGSS]  G.Baumslag, S.Gersten, M.Shapiro and H.Short, *Automatic groups and amalgams*, J. of Pure and Appl. Algebra **76** (1991), 229-316.

[Br]    M.R.Bridson, *Geodesics and curvature in metric simplicial complexes*, Group theory from a geometric viewpoint, Proc. ICTP. Trieste, World Scientific, Singapore, 1991.

[BH]    M.R.Bridson and A.Haeflinger, *An introduction to CAT(0)-spaces*, in preparation.

[Com1]  L.Comerford, Jr., *Powers and conjugacy in small cancellation groups*, Arch. Math. (Basel) **26** (1975), no. 4, 353–360.

[Com2]  L.Comerford, Jr., *A note on power-conjugacy*, Houston J. Math. **3** (1977), no. 3, 337–341.

[C1]    G.Conner, *Metrics on Groups*, PhD Thesis, University of Utah, 1992.

[C2]    G.Conner, *Properties of translation numbers in solvable groups*, Brigham Young University, preprint (1994).

[C3]    G.Conner, *A finitely generated group with irrational translation numbers*, Brigham Young University, preprint (1994).

[CT]    L.Comerford, Jr. and B.Truffault, *The conjugacy problem for free products of sixth-groups with cyclic amalgamation*, Math. Z. **149** (1975), no. 2, 169-181.

[ECHLPT] D.B.A.Epstein, J.W.Cannon, D.F.Holt, S.V.F.Levy, M.S.Paterson and W.P.Thurston, *Word Processing in Groups*, Jones and Bartlett, MA, 1992.
[GS1] S.Gersten and H.Short, *Rational Subgroups of Biautomatic Groups*, Ann. Math. **134** (1991), 125-158.
[GS2] S.Gersten and H.Short, *Small cancellation theory and automatic groups*, Invent. Math. **102** (1990), 305-334.
[GS3] S.Gersten and H.Short, *Small cancellation theory and automatic groups:Part II*, Invent. Math **105** (1991), 641-662.
[Gr] M.Gromov, *Hyperbolic Groups*, in 'Essays in group theory', edited by S.M.Gersten, MSRI Publ. 8, Springer, 1987.
[J1] K. Johnsgard, *Geodesic tilings for equal geodesic words in $C''(p) - T(q)$ group presentations*, Cornell University, preprint, 1993.
[J2] K. Johnsgard, *Two automatic spanning trees in small cancellation group presentations*, Intern. J. Alg. Comp. (to appear).
[L] S.Lipschutz, *On Greendlinger Groups*, Comm. on Pure and Appl. Math **XXIII** (1970), 743-747.
[LS] R.Lyndon and P.Schupp, *Combinatorial Group Theory*, Springer-Verlag, New York, 1977.
[MP] M.El-Mosalamy and S.Pride, *On $T(6)$ groups*, Math. Proc. Camb. Philos. Soc. **102** (1987), 443-451.
[S] E.Swensen, *Hyperbolic Elements in Negatively Curved Groups*, Geom. Dedicata **55** (1995), 199–210.
[W] C.Weinbaum, *The word and conjugacy problem for the knot group of any prime alternating knot*, Proc. Am. Math. Soc. **22** (1971), 22-26.

DEPARTMENT OF MATHEMATICS, GRADUATE SCHOOL AND UNIVERSITY CENTER OF THE CITY UNIVERSITY OF NEW YORK, 33 WEST 42-ND STREET, NEW YORK, NY10036
*Current address*: Department of Mathematics, City College of CUNY, Convent Avenue at 138-th Streer, New York, NY10031
*E-mail address*: ilya@groups.sci.ccny.cuny.edu